\input amstex
\documentstyle {amsppt}
%\magnification=1200
\tolerance=3000
\openup 6pt
\nologo

\topmatter
\title
The existence of $\sigma-$finite invariant measures,\\
Applications to real 1-dimensional dynamics\\
\endtitle

\author
Marco Martens
\endauthor

\affil
IMPA, Estrada Dona Castorina 110, Rio de Janeiro, Brasil 
\endaffil

%\date{August 1991}

\abstract{ A general construction for  $\sigma-$finite absolutely continuous
invariant measure will be presented. It will be shown that the local bounded
distortion of the Radon-Nykodym derivatives of $f^n_*(\lambda)$ will imply
the existence of a $\sigma-$finite invariant measure for the map $f$ 
which is absolutely
continuous with respect to $\lambda$, a measure on the phase space describing
the sets of measure zero.

\flushpar
Furthermore we will discuss sufficient conditions for the existence of
$\sigma-$finite invariant absolutely continuous measures for real 1-dimensional
dynamical systems.}
\endabstract

\endtopmatter

\openup 6pt
\bigskip
\centerline{\bf 1. Introduction}
\bigskip
The statistical study of a dynamical system begins with the question whether or
not the system has an absolutely continuous invariant measure, finite or
$\sigma-$finite. In 1947 Halmos gave a characterization
 of the (bijective)
dynamical systems which have a $\sigma-$finite absolutely continuous invariant
measure, see [Ha]. During this time there was some hope that every dynamical system has
$\sigma-$finite invariant measures. Unfortunately this turned out not to be 
true. Ornstein gave an
example of a piecewise linear bijective map on the interval not having such a
measure ([O]). 

\flushpar
Here we will give a characterization of the ergodic conservative dynamical
systems on locally compact spaces having $\sigma-$finite absolutely continuous
invariant measures.       

\flushpar
The origin of the characterization presented here can be found in the theory of
real 1-dimensional dynamics and the theory of Markov processes. In [H] Harris
used limits of ratios of long term transition probabilities to construct
infinite stationary states for Markov processes on countable state spaces. In
section 2 we will use this idea for constructing $\sigma-$finite
invariant measures. The construction gives rise to the existence
theorem A. The distortion of a measure, used in the theorem, will be defined
precisely in section 2. Furthermore remember that a map $f:X\to X$ is ergodic 
conservative with respect to a measure $\lambda$ on $X$ if every set of
positive measure is hit by the orbits of $\lambda-$almost all points (see [P]).

\proclaim{Theorem A} Let $\lambda$ be a Borel probability measure on the 
$\sigma-$compact space $X$. The ergodic conservative map $f:X \to X$  has a 
$\sigma-$finite invariant measure absolutely
continuous with respect to $\lambda$ if the
Radon-Nykodym derivatives of $f^n_*\lambda$ have locally bounded distortion.
\endproclaim

\flushpar
As in the general construction of invariant probability measures the
construction is done by pushing forward some initial measure and then
considering limits of these push-forwards. It turns out that the procedure only
gives rise to $\sigma-$finite absolutely continuous invariant measures if the
initial measures are of some special type.  In section 3  we will construct the
initial measures.
 
\flushpar
In section 4  we will study the existence of $\sigma-$finite invariant measures
for real 1-dimensional differentiable dynamics. As we know from [J] there is 
no general
existence theorem  for absolutely continuous invariant probability measures:
 there exist conservative  quadratic maps on the interval not
having
absolutely continuous invariant probability measures. Even the existence
question for $\sigma-$finite  absolutely continuous invariant  measures
can not be
answered in general.  Katznelson ([Ka]) constructed diffeomorphisms of the
circle not having $\sigma-$finite absolutely continuous invariant measures. 

\flushpar
Applying the developed theory we can formulate sufficient
conditions implying the existence of $\sigma-$finite absolutely continuous
invariant measures for real 1-dimensional differentiable dynamics. 

\flushpar
In [HKe2] Hofbauer and Keller gave an existence theorem for some type of conservative
unimodal maps. Now this theorem  can be generalized to multimodal and also
dissipative maps:

\proclaim{Theorem B}  Let $f$ be a $C^3$ map on the interval (or the circle) 
satisfying
\parindent=15pt
\item{1)} $f$ has only finitely many critical points, points where the 
derivative vanishes, and the Schwarzian derivative is everywhere negative
except in the critical points;
\item{2)} there exists a dense orbit;
\item{3)} the orbits of the critical points stay in a closed invariant 
  set of Lebesgue
measure zero.

\flushpar
Then $f$ has a $\sigma-$finite absolutely continuous invariant measure.
\endproclaim 

\flushpar
In [HKe1] quadratic maps are shown to exist having very strange Bowen-Ruelle
measures. The same techniques can be used to show that there exists a quadratic 
map
whose critical orbit is in a Cantor set but which doesn't have an absolutely
continuous invariant  probability measure. This means that in general the
invariant measures of Theorem B are really $\sigma-$finite. Furthermore we 
obtain a $\sigma-$finite Folklore theorem.  

\flushpar
All the results
concerns maps whose critical orbits stay in some closed invariant 
 set of Lebesgue measure
zero. In the other case, some critical orbits are dense, we  state the 

\proclaim{Conjecture} There exist conservative  quadratic maps on the interval not having
$\sigma-$finite absolutely continuous invariant measures.
\endproclaim  

\flushpar
In the appendix we will give  a short proof of the Chacon-Ornstein Theorem, 
the main theorem in $\sigma-$finite ergodic theory.

\flushpar
Remember the following notation: if $g:X\to X$ is a Borel measurable function 
and $\mu$ a Borel measure on $X$ then $g_*\mu$ is the measure defined by
$g_*\mu(A)=\mu(g^{-1}(A))$.

%\magnification=1200
\tolerance=3000

\bigskip
\centerline{\bf 2. The construction of $\sigma-$finite absolutely continuous
invariant measures}
\bigskip
In this section we are going to  construct
 $\sigma-$finite absolutely continuous invariant measures. Let $\lambda$ be 
the Borel measure describing the sets of measure zero. A 
starting point for constructing absolutely continuous invariant measures is
considering limits of the Birkhoff sums $\{\frac1n\Sigma_{i=0}^{n-
1}f_*^i\lambda\}_{n\ge 0}$. Indeed, if for all sets $A$ with $\lambda(A)>0$ the
measures $f^n_*\lambda(A)$ stay away from 0 we can construct an absolutely
continuous invariant probability measure, simply by taking a converging
subsequence of Birkhoff sums. 
In case our system doesn't have an
invariant probability measure, that is there exist closed sets $A$ with 
$f^n_*\lambda(A)$
converging to zero, we have to consider other limits. 

\flushpar
In [H] Harris used limits of ratios of long term transition probabilities
to construct stationary states for Markov processes. In our construction we
will choose some set $I_0$ of positive measure and consider limits of 
the normalized sequence 
$$
\frac{\sum_{i=0}^{n-1}f^i_*\lambda}
     {\sum_{i=0}^{n-1}f^i_*\lambda(I_0)}.     
$$
Our construction is strongly related to the Chacon-Ornstein Theorem or, which
is in some sense equivalent to it, the Doeblin-Ratio-Limit Theorem for Markov
processes (see resp. [K],[P] and [F]). 

\bigskip
\flushpar
First let us remember some general notions. Let $X$ be a $\sigma-$compact
topological space, it can be written as a countable union of compact sets, and
$\Cal{B}(X)$ be the set of Borel measures on $X$. The set $\Cal{B}_\sigma(X)$
consists
of all measures $\mu \in \Cal{B}(X)$ for which there exists a collection
$\{X_n|n \in \bold{N}\}$ of pairwise disjoint measurable sets in $X$ such
that
$$
\align
&\mu(X_n)\text{ is finite for all } n \in \bold{N};\\
&\mu(X-(\bigcup_{n\in \bold{N}} X_n))=0.
\endalign
$$
The measures in $\Cal{B}_\sigma(X)$ are called $\sigma-${\it finite
measures\/} on $X$.

\bigskip
\flushpar
For discussing the notion of absolute continuity we fix a measure $\lambda \in
\Cal{B}_\sigma(X)$ determining the {\it null sets\/}, the sets which have to
have
measure zero. We may assume without restricting generality that $\lambda$ is a
probability measure. Remember a measure $\mu \in \Cal{B}_\sigma(X)$ is called
absolutely
continuous with respect to $\lambda$, $\lambda \gg \mu$, iff the sets of
$\lambda-$measure zero also have $\mu-$measure zero. These absolutely
continuous measures can be expressed by an integral: for all $\mu$ with 
$\lambda \gg \mu$
there exists an essentially  unique integrable non-negative function $\rho$ such that
$\mu(A)= \int_A \rho
d\lambda$ for all measurable sets $A\subset X$. This function is called the
{\it Radon-Nykodym derivative} (or just derivative) of $\mu$ with respect to
$\lambda$. 

\flushpar
The main objects to be studied here will be absolutely continuous measures
whose Radon-Nykodym  derivatives has locally bounded distortion: we say that
the derivative of $\mu \in \Cal{B}_\sigma(X)$, $\lambda \gg \mu$, on $I\subset
X$ with $\mu(I)$ finite,
has {\it distortion} bounded by $K$ iff for all measurable sets $A \subset I$ 
$$
\frac1K\frac{\lambda(A)}{\lambda(I)}\le \frac{\mu(A)}{\mu(I)}
                             \le K \frac{\lambda(A)}{\lambda(I)}.
$$
Observe that the constant $K$ can be taken to be equal to $1$ iff $\mu$
has constant derivative on $I$. If you consider a Radon-Nykodym derivative
as
an object deforming the original measure $\lambda$ this notion of distortion
is related to the concept of distortion for differentiable maps of the
interval (see [GuJ], [LB], [MMS] or [Sw]).

\bigskip
\flushpar
In this section we fix a probability measure $\lambda \in \Cal{B}_\sigma(X)$
and a measurable map
$f:X\rightarrow X$ on the $\sigma-$compact space $X$. Assume that 
$\lambda$ is quasi-invariant for $f$:
$\lambda \gg f_*\lambda$.  

\bigskip
\flushpar
The next step is to give some definitions which enables us to deal with
$\sigma-$finite absolutely continuous invariant measures, shortly {\it acim
\/}
(if we want to emphasize that some acim is a probability measure 
we call it {\it acip\/}: absolutely continuous invariant probability measure).

\flushpar
A $\lambda-${\it partition\/} $\Cal{G}$ of $X$ is a countable collection of
pairwise disjoint Borel sets of $X$, say $\Cal{G}=\{I_n|n\in \bold{N}\}$,
such
that for all $n \in \bold N$
\parindent=15pt
\item{1)} $I_n$ is $\sigma-$compact;
\item{2)} $0<\lambda(I_n)<\infty$;
\item{3)} $\lambda(X-(\bigcup_{I\in \Cal{G} }I))=0$.

\flushpar
If the $\lambda-$partition $\Cal{G}$ of $X$ has the additional property
\parindent=15pt
\item{4)} for all pairs $I_1,I_2 \in \Cal{G}$ there exists $n \ge 0$ such that
$\lambda(f^{-n}(I_1)\cap I_2)>0$

\flushpar
we will say $f$ is $\Cal{G}-$irreducible.

\flushpar
The role of the $\lambda-$partition is the following: its elements will turn
out to be sets of finite measure for the acim of the $\Cal{G}-$irreducible
 map $f$.

\bigskip
\flushpar
In the sequel we fix a $\lambda-$partition $\Cal{G}$ and assume that
$f:X\rightarrow X$ is $\Cal{G}-$irreducible. 

\bigskip
\flushpar
To define the measure spaces
$\Cal{M}(\Cal{G},f), \Cal{M}_s(\Cal{G},f)$ and $\Cal{M}_{\infty}(\Cal{G},f)$,
in which the construction of the acims will take place, we need some
properties of measures $\mu \in \Cal{B}_\sigma(X)$.

\flushpar
Let $I\in \Cal{G}$ and $K>0$:
\parindent=50pt
\item{m1($I,K$):} $\mu(I) \in [\frac1K,K]$;  
\item{m2($I,K$):} For all $n \ge 0$ the $\mu-$measure of $f^{-n}(I)$ is finite
and positive and the derivatives of the measures $f^{n}_*\mu$ on $I$ with 
respect to $\lambda$ has
distortions bounded by $K$.  This property states that measures $f^n_*\mu$
have locally uniformly bounded  distortions.
\item{m3($I$):} $\sup_{n\ge0}\mu(f^{-n}(I)) < \infty$;
\item{m4($I$):} $\dsize\sum_{n=0}^{\infty}\mu(f^{-n}(I))=\infty$.

\bigskip
\flushpar
The measure spaces we need are
$$
\align
\Cal{M}(\Cal{G},f)&=  \{\mu \in \Cal{B}_\sigma(X)| \text{ for all $I\in
\Cal{G}$ 
there is a $K>0$ with m1($I,K$) and m2($I,K$) }\};\\
\Cal{M}_s(\Cal{G},f)&=\{\mu \in \Cal{M}(G,f)| \text{ m3($I$) holds for all  } 
I \in \Cal{G}\};\\
\Cal{M}_{\infty}(\Cal{G},f)&=\{\mu \in\Cal{M}_s(\Cal{G},f)| \text{ m4($I_0$)
holds for some } I_0 \in \Cal{G}\}.
\endalign
$$

\flushpar
Because we are only considering the fixed $\lambda-$partition $\Cal{G}$ and
the fixed $\Cal{G}-$irreducible map $f:X\rightarrow X$ we will use the short 
names
$\Cal{M}$, $\Cal{M}_s$ and $\Cal{M}_{\infty}$. Furthermore observe that
$\Cal{M}_{\infty} \subset \Cal{M}_s\subset \Cal{M}$ and that 
$\Cal{M}$ only contains measures which are absolutely 
continuous with respect to $\lambda$.

\bigskip
\flushpar
The construction of the acims will be by taking converging subsequences in
$\Cal{M}$.
Hence we have to describe which kind of convergence  we are going to use. A
sequence $\mu_n$ in $\Cal{M}=\Cal{M}(\Cal{G},f)$ is said to converge to $\mu
\in \Cal{B}_\sigma(X)$ iff for all $I \in \Cal{G}$ and for every compact
$A\subset I$ 
we have weak convergence of $\mu_n|A \rightarrow \mu|A$.

\flushpar
For two reasons the space $\Cal{M}$ is not compact: the measures are not
assumed to be bounded and the underlying space is not compact; mass can
disappear to the boundary of the space. Let us try to describe some
compact subsets of $\Cal{M}$. A collection $\Cal{A} \subset \Cal M$ is called 
{\it uniform} iff for all $I\in \Cal G$ there exists $K(I)>0$  such that
m1($I,K(I)$) and m2($I,K(I)$) hold for all $\mu \in \Cal A$.

\flushpar
The campactness in the weak topology of the set of probability measures on a
compact space implies easily

\proclaim{Lemma 2.1} Uniform collections in $\Cal M(\Cal G,f)$ have compact
closures in $\Cal M(\Cal G,f)$.
\endproclaim

\flushpar
Now we are going to use the above measure spaces to construct the
acims.
Let $\mu \in \Cal M$ and define the following measures in $\Cal B_\sigma(X)$:
$$
\align
&S_n\mu=\sum_{i=0}^{n-1}f^i_*\mu;\\
&Q_n\mu=\frac{S_n\mu}{S_n\mu(I_0)}.
\endalign
$$
for all $n \ge 0$. $I_0$ is a fixed element of $\Cal{G}$. If $\mu \in \Cal
M_{\infty}$ then we choose $I_0$ such that m4($I_0$) holds.

\flushpar
The lemma which will assure the existence of limits is

\proclaim{Lemma 2.2} Let $\mu \in  \Cal M_s(\Cal G,f)$ then the collection
$$
\Cal A_{\mu}=\{Q_n\mu|n \ge 0\}
$$
is uniform.
\endproclaim

\flushpar
Before proving this lemma we are going to use lemma 2.1 and 2.2 to define the
following limit set in $\Cal{M}$. Let $\mu \in \Cal{M}_s$ then
$$
\omega(\mu) \subset \Cal{M}(\Cal{G},f)
$$
is the set of all limits of the sequence $\Cal{A}_{\mu}$. We will look for
acims in these sets $\omega(\mu)$. For the moment we know already
that it only contains measures which are equivalent to $\lambda$.

\flushpar
The next lemma shows how the measures of backward orbits of two elements in
$\Cal{G}$ are related. It will be used at several places; it serves for
gluing
together the information given by the local boundedness of the distortions.  

\proclaim{Lemma 2.3} Let $\mu \in \Cal{M}(\Cal{G},f)$. Then for every pair
$I_1,I_2 \in \Cal G$ there exist $\epsilon >0$ and $n_0 \ge 0$ such that
$$
\mu(f^{-n-n_0}(I_1))\ge \epsilon \mu(f^{-n}(I_2)) 
$$
for $n\ge 0$.
\endproclaim
\demo{proof} Because $f$ is $\Cal G-$irreducible there exists $n_0\ge 0$ such that
$\lambda(f^{-n_0}(I_1)\cap I_2) > 0$. Hence for $\mu \in \Cal{M}(\Cal{G},f)$
we
get for all $n\ge 0$ 
$$
\align
\mu(f^{-n-n_0}(I_1))&\ge\frac{\mu(f^{-n}(f^{-n_0}(I_1)\cap I_2))}
                             {\mu(f^{-n}(I_2))}
                        \mu(f^{-n}(I_2))\\
                    &\ge\frac1K\frac{\lambda(f^{-n_0}(I_1)\cap I_2)}
                                    {\lambda(I_2)}
                                \mu(f^{-n}(I_2))\\
                   &=\epsilon \mu(f^{-n}(I_2)),
\endalign
$$
where $K>0$ is such that m2($I_2,K$) holds for $\mu$.
\qed
\enddemo

\flushpar
This lemma is the place where we use the $\Cal G$-irreducibility of $f$.

\proclaim{Lemma 2.4} Let $\mu \in \Cal{M}_s(\Cal{G},f)$. For every pair
$I_1,I_2 \in \Cal{G}$ there exists $K<\infty$ such that
$$
\frac1K \le \frac{S_n\mu(I_1)}{S_n\mu(I_2)}\le K
$$
for all $n\ge 0$.
\endproclaim
\demo{proof}
Let $\epsilon, n_0$ be given by lemma 2.3. For $n\le n_0$ we have some bound.
Let $n>n_0$. Then
$$
\align
\frac{S_n\mu(I_1)}{S_n\mu(I_2)}
&\ge \frac{\sum_{i=0}^{n-n_0-1}\mu(f^{-i-n_0}(I_1))}
          {\sum_{i=n-n_0}^{n-1}\mu(f^{-i}(I_2))+
           \sum_{i=0}^{n-n_0-1}\mu(f^{-i}(I_2))}\\
&\ge \epsilon
     \frac{\sum_{i=0}^{n-n_0-1}\mu(f^{-i}(I_2))}
          {n_0 \sup_{i \ge 0}\mu(f^{-i}(I_2))+
           \sum_{i=0}^{n-n_0-1}\mu(f^{-i}(I_2))}\\
&\ge \epsilon
     \frac{\mu(I_2)}
          {n_0 \sup_{i\ge 0}\mu(f^{-i}(I_2)) +\mu(I_2)}
\endalign
$$
which is a finite positive number. Remark that we used in the last step that
the function $x\to \frac{x}{a+x}$ is increasing.

\flushpar
By interchanging the role of $I_1$ and $I_2$ we also get an upper bound.
\qed
\enddemo

\demo{proof of lemma 2.2} Fix $I\in \Cal{G}$.
\demo{m1} Let $K$ be the number given by lemma 2.4 applied to $I$ and $I_0$.
We
get directly from lemma 2.4 that m1($I,K$) holds for all $Q_n\mu$, $n\ge 0$.
\enddemo
\demo{m2} Let $K>0$ be such that m2($I,K$) holds for $\mu$. Fix $m\ge 0$ and
let $A\subset I$ be measurable. Then we get for $n\ge 0$
$$
\align
\frac{Q_m\mu(f^{-n}(A))}{Q_m\mu(f^{-n}(I))}
&=\frac{S_m\mu(f^{-n}(A))}{S_m\mu(f^{-n}(I))}\\
&=\frac{\sum_{i=n}^{m+n-1}\mu(f^{-i}(A))}
       {\sum_{i=n}^{m+n-1}\mu(f^{-i}(I))}.
\endalign
$$
Because m2($I,K$) holds for $\mu$ we easily get that this last number is in
the
interval $[\frac1K\frac{\lambda(A)}{\lambda(I)},
                 K \frac{\lambda(A)}{\lambda(I)}]$. This proves the lemma.
\qed
\enddemo

\flushpar
The  following lemma tells under which conditions on $\mu$, $\omega(\mu)$ will
contain
invariant measures.

\proclaim{Lemma 2.5} Let $\mu \in \Cal{M}_{\infty}(\Cal{G},f)$. Then
$\omega(\mu)$ contains only invariant measures.
\endproclaim
\demo{proof}
Let $\nu \in \omega(\mu)$, say $\nu=\lim Q_n\mu$ ($\lim$ means : the limit of
a
certain converging subsequence). Because $\lambda(X-\cup\Cal{G})=\nu(X-
\cup\Cal{G})=0$ and $\lambda$ is quasi-invariant for $f$ we only have
to consider $A\subset I$, $I\in \Cal{G}$. Let $A\subset I$ be compact. Then
$$
\align
\nu(f^{-1}(A))
&=\lim \frac{S_n\mu(f^{-1}(A))}
            {S_n\mu(I_0)}\\
&=\lim \frac{\sum_{i=0}^{n-1}\mu(f^{-i}(A))-
             \mu(A)+\mu(f^{-n}(A))}
            {\sum_{i=0}^{n-1}\mu(f^{-n}(I_0))}\\
&=\nu(A)+\lim\frac{\mu(f^{-n}(A))-\mu(A)}
                  {S_n\mu(I_0)}.
\endalign
$$
Now we use that m3($I$) and m4($I$) hold for $\mu$ and we get
$$
\nu(f^{-1}(A))=\nu(A). 
$$
The measure is invariant.
\qed
\enddemo 

\proclaim{Proposition 2.6} Let $\lambda \in \Cal{B}_\sigma(X)$ be a Borel
measure on
the
$\sigma-$compact space $X$ and $f:X\rightarrow X$ a measurable map. The map $f$ 
has a $\lambda-$equivalent $\sigma-$finite invariant measure if 
it is $\Cal{G}-$irreducible for some $\lambda-$partition $\Cal{G}$ of $X$ with
$$
\Cal{M}_{\infty}(\Cal{G},f)\ne \emptyset.   
$$ 
The elements of $\Cal{G}$ will be pieces of $X$ with bounded measure.
\endproclaim

\demo{proof} If $\mu \in\Cal{M}_{\infty}(\Cal{G},f)$ then from lemma 2.2, 2.1
 we get $\omega(\mu)\ne \emptyset$. Furthermore Lemma 2.5 tells us that 
$\omega(\mu)$ only contains invariant measures.
\qed
\enddemo

\bigskip
\flushpar
In fact we also want the reverse statement: $f$ has an acim iff
it is $\Cal{G}-$irreducible for some $\lambda-$partition $\Cal{G}$ of $X$ with
$\Cal{M}_{\infty}(\Cal{G},f)\ne \emptyset$.

\flushpar
It is not hard to prove that for every map $f:X\rightarrow X$ which has an
acim there exists a $\lambda-$partition $\Cal{G}$ such that
$\Cal{M}_{\infty}(\Cal{G},f)\ne \emptyset$. A
problem arises when we want to get it such that $f$ becomes 
$\Cal{G}-$irreducible.
Probably it is possible to get this property. This technical problem can be
illustrated by the question:  does the feigenbaum map have an acim?
\flushpar
We can overcome this technical problem by assuming that $f$ is ergodic and
conservative with respect to $\lambda$: every set of positive
$\lambda-$measure will intersect every other set of positive $\lambda-$measure
after some time. If $f$ is ergodic and conservative it will be
$\Cal{G}-$irreducible for every $\lambda-$partition $\Cal{G}$.

\flushpar
So we get the following: an ergodic conservative map $f$ has an acim iff there exists a
$\lambda-$partition $\Cal{G}$ such that $\Cal{M}_{\infty}(\Cal{G},f)\ne
\emptyset$.

\flushpar
Using the following lemma we even can state a stronger existence theorem.

\proclaim{Lemma 2.7} Let $\mu \in \Cal{B}_\sigma(X)$ with $\lambda \gg \mu$. If
$f$ is ergodic  and conservative with respect to
$\lambda$ then for every set $A$ with $\lambda(A)>0$
$$
\sum_{i=0}^\infty \mu(f^{-i}(A))=\infty.
$$ 
\endproclaim
 
\demo{proof}
The ergodicity and conservatively tells us that almost every point in $A$ will
return to $A$ infinitely many times. Now use the Borel-Cantelli Lemma.
\qed
\enddemo

\proclaim{Corollary 2.8} An ergodic conservative map $f$ has a $\sigma-
$finite absolutely continuous invariant measure iff there exists a partition
$\Cal{G}$ with
$$
\Cal{M}_s(\Cal{G},f)\ne \emptyset.
$$
\endproclaim

\flushpar
The precise formulation of Theorem A in the introduction goes as follows.

\proclaim{Theorem 2.9} An ergodic conservative map $f$ has a $\sigma-$finite
absolutely continuous measure if there exists a $\lambda-$partition $\Cal{G}$  with
$$
\lambda \in \Cal{M}(\Cal{G},f).
$$
\endproclaim

\flushpar
This means: once the derivatives of $f^n_*\lambda$ have locally bounded
distortion the existence of an acim is assured.

\bigskip
\centerline{\bf 3. The initial measures}
\bigskip

\parindent=30pt 
In this section we are going to study a condition which will imply 
$\Cal{M}_\infty(\Cal{G},f)\ne \emptyset$, that is, it will imply the existence
of acims.

\flushpar
Fix $\lambda \in \Cal{B}_\sigma(X)$ and a $\lambda-$partition $\Cal{G}$ of
$X$.
We say
that the measurable map $f:X\rightarrow X$ is {\it finite-to-1 with respect to
} $\Cal{G}$ iff for all $I\in \Cal{G}$ $f^{-1}(I)$ is up to a nullset 
contained in a finite subcollection of $\Cal{G}$. Furthermore
$PL_{\lambda}(\Cal{G})$ is the set of all {\it distortion free } measures:
these measures are equivalent to $\lambda$ with densities which are constant
on
the element of $\Cal{G}$.

\proclaim{proposition 3.1} Let $\lambda \in \Cal{B}_\sigma(X)$ be a Borel
measure
on $X$
and $\Cal{G}$ a $\lambda-$partition such that the measurable map
$f:X\rightarrow X$ is $\Cal{G}-$irreducible.

\flushpar
If there exists a $\lambda-$partition $\Cal{G}_0$ with
\parindent=15pt
\item{1)} $\Cal{G}$ is a refinement of $\Cal{G}_0$;
\item{2)} $PL_{\lambda}(\Cal{G}_0)\subset \Cal{M}(\Cal{G},f)$ 

\flushpar
then
$$
\Cal{M}_{\infty}(\Cal{G},f)\ne \emptyset.
$$
\endproclaim

\demo{proof}
The condition $PL_{\lambda}(\Cal{G}_0)\subset \Cal{M}(\Cal{G},f)$ is a strong
condition. It implies that $f$ is finite-to-1 with respect to $\Cal{G}_0$: if
$f$
is not finite-to-1 it is easy to find $I\in \Cal{G}_0$ and $\mu \in
PL_{\lambda}(\Cal{G}_0)$ such that $\mu(f^{-1}(I))=\infty$. 

\flushpar
We are going to define a measure $\mu \in PL_{\lambda}(\Cal{G}_0)$ satisfying
m3($I_0$) and m4($I_0$), where $I_0 \in \Cal{G}$ is fixed. Using
$PL_{\lambda}(\Cal{G}_0) \subset \Cal{M}(\Cal{G},f)$ and lemma 2.4 we get $\mu
\in \Cal{M}_\infty(\Cal{G},f)$.

\flushpar
Let $\Cal{G}_0=\{J_n|n\ge 0\}$ and $\Cal G=\{I_n|n\ge 0\}$. We can assume $I_0\subset J_0$. Define
$L_N=\cup_{i=0}^N J_i$ for $N\ge 0$. We are
going to define $\mu$ by giving its density $\delta$ with respect to $\lambda$
$$
\delta= \sum_{N \ge 0} c_N \bold{1}_{J_N}. 
$$ 
The numbers $c_N > 0$ will be defined inductively satisfying the following 
induction hypothesis
$$
\sup_{n\ge 0} \mu|L_N(f^{-n}(I_0))=1-(\frac12)^{N+1}.
$$ 
Because $\lambda(I_0)<\infty$ we can choose $c_0>0$ such that the induction
hypothesis holds for $N=0$. Suppose that $c_0,c_1,\dots, c_N$ are defined
satisfying the induction hypothesis. This means that the measure $\mu|L_N$ is
well defined. Now we have to define the value $c_{N+1}$ of the density on
$J_{N+1}$: let
the map $c \to \mu_{N,c} \in PL_{\lambda}(\Cal{G}_0)$ with $c \in [0,\infty)$
be
defined as follows
$$
\mu_{N,c}=\mu|L_N+ c\lambda|J_{N+1}.
$$
Using the fact that $f$ is $\Cal{G}-$irreducible, that is there exists an $n\ge 0$
such that $\lambda(f^{-n}(I_0) \cap J_{N+1})\ne 0$, it is easy to see that the
map $\phi:[0,\infty) \to \bold{R}$ defined by
$$
\phi(c)=\sup_{n\ge 0} \mu_{N,c}(f^{-n}(I_0))
$$
tends continuously to infinity for $c\to \infty$. From the definition of
$c_0,c_1,\dots, c_N$ we get $\phi(0)=1-(\frac12)^{N+1}$. Hence there exists
$c_{N+1}>0$ such that 
$$
\sup_{n\ge 0} \mu|L_{N+1}(f^{-n}(I_0))= \phi(c_{N+1})=1-(\frac12)^{N+2}.
$$
We finished the induction step; the measure $\mu$ is well defined.

\flushpar
\demo{proof of m3($I_0$)}
Suppose there exists $n\ge 0$ such that $\mu(f^{-n}(I_0))>1$. Because
$\Cal{G}_0$
is an exhausting partition of $X$ there exists $N\ge 0$ such that
$$
\mu|L_N(f^{-n}(I_0))>1
$$ 
which contradicts the construction of $\mu$.
\enddemo
\demo{proof of m4($I_0$)}
We are going to construct a sequence $n_k \to \infty$ such that
$$
\mu(f^{-n_k}(I_0))\ge \frac12.
$$
This will imply m4($I_0$).

\flushpar
Suppose we have a finite set $\{n_i|i=0,1,\dots,k\}$ such that for all of them

$\mu(f^{-n_i}(I_0))\ge \frac12$. Let us find another one having this property.
Because $f$ is finite-to-1 with respect to $\Cal{G}_0$ there exists $N\ge 1$
such
that 
$$
\bigcup_{i=0}^k f^{-n_i}(I_0) \subset L_N   
$$
up to a set of measure zero.
Hence $\mu(f^{-n_i}(I_0))\le 1-(\frac12)^{N+1}$ for  $i=0,1,\dots,k$. From the
definition of $c_{N+1}$ we easily get a number $n_{k+1}$ such that
$$
1-(\frac12)^{N+1} < \mu(f^{n_{k+1}}(I_0)) 
$$
which is obviously not one of the previous ones.
\qed
\enddemo
\enddemo

\flushpar
Observe that theorem 2.9 gives a much weaker sufficient condition for the
existence of acims for conservative maps. The use of proposition 3.1 will be
for general maps. Indeed it can be shown that lemma 2.7 gives a
characterization for dissipative unimodal maps: A unimodal map is dissipative
iff $\Sigma_{i=0}^{\infty}\lambda(f^{-i}(A))<\infty$ for all sets $A \subset
\cup \Cal{G}$ where $\Cal{G}$ is some $\lambda-$partition.
%\magnification=1200
\tolerance=3000

\bigskip
\centerline{\bf 4. Applications to 1-dimensional real dynamics}
\bigskip

In this section we will discuss the existence of absolutely continuous
invariant measures for maps on the interval having negative Schwarzian
derivative. The existence of invariant probability measures is strongly related
to the expansion along the orbits of the critical points, the points where the
derivative vanishes. In [CE] the existence of acips was shown for unimodal maps
having exponential
growth of the derivative  along the critical orbit. In [NS] this result was
obtained for a weaker growth of the derivative, for example (non-linear) 
polynomial growth turns
out to be sufficient.

\flushpar
Another type of existence theorems is described in [Mi] and [S]: if the
orbits of the critical points are not accumulating at critical points, maps
having this property are called {\it Misiurewicz} maps, then an acip exists. 
Here we will describe some results continuing in this direction.

\flushpar
In general the orbit closures of the critical points of Misiurewicz maps will lie in
some closed invariant set of measure zero which doesn't contain critical
points.
In the sequel we will allow our maps to exhibit some recurrence, critical
orbits accumulate at critical points but we continue imposing the critical
orbits  to be in some closed invariant set of measure zero. For this type of
maps the existence of acims will be shown. A result in this direction was
already obtained in [HKe2]: if the critical point of a conservative unimodal map
with negative Schwarzian derivative stays in a Cantor set  then an acim
exists.  

\flushpar
The question whether these unimodal maps are always conservative is the main
open problem in the theory of interval dynamics. A natural candidate for
having an absorbing Cantor set was recently proved to be
conservative, see [LM]. 

\flushpar
The results presented in this section can be formulated shortly as follows:
multimodal maps whose critical orbits lie in a closed  invariant set have a
$\sigma-$finite absolutely continuous invariant measure. In particular,
unimodal maps having an absorbing Cantor set have an acim.  

\bigskip
\flushpar
In the sequel $X$ denotes the interval $[0,1]$ or the circle endowed with the
Lebesgue measure $\lambda$. 

\bigskip
\flushpar
Let us first define  the main analytical tool we need. Let $g:I\to J$ be $C^3$
and  mapping the interval $I\subset X$ to the interval $J\subset X$. The Schwarzian derivative $Sg:I\to \bold{R}$ of $g$ is defined to be
$$
Sg(x)=\frac{D^3g(x)}{Dg(x)}-\frac32(\frac{D^2g(x)}{Dg(x)})^2.
$$
\flushpar
An important and easy to derive property of maps with negative Schwarzian
derivative is that the iterates of these maps also have negative Schwarzian derivative.

The Schwarzian derivative enables us to formulate the following distortion result.

\proclaim{Koebe-Lemma} For every $\epsilon>0$ there exists $K>0$ with the following property. Let $g:I\to J$ be a diffeomorphism  mapping the interval $I\subset X$ to the interval $J\subset X$. Assume that $Sg(x)<0$ for all $x\in I$.

\flushpar
If  $M\subset I$ is an interval such that the components of $I-M$, denoted by
$L$ and $R$, satisfy
$$
\frac{\lambda(g(L))}{\lambda(g(M))}\ge \epsilon
\text{ and }
\frac{\lambda(g(R))}{\lambda(g(M))}\ge \epsilon
$$
then
$$
\frac1K\le \frac{|Dg(x_1)|}{|Dg(x_2)|}\le K
$$
for all $x_1,x_2\in M$.
\endproclaim

\flushpar
The proof of this fundamental lemma can be found in different places ([GuJ], [MMS]).

\bigskip
\flushpar
The class $\Cal D(X)$ of functions which we are going to consider is the class
of piecewise diffeomorphic maps on $X$.  These maps are defined as follows. Let
$f\in \Cal D(X)$ then there exists a $\lambda-$partition $\Cal P$ of $X$
consisting of open intervals such that for all $I\in \Cal P$ the restriction 
$f|I$ is a diffeomorphism with negative Schwarzian derivative. Furthermore we
assume $f$ to have a dense orbit. 

\flushpar
An open interval $T$ is called a branch of $f^n$, $f\in \Cal D(X)$, if $T$ is a
maximal interval on which $f^n$ is diffeomorphic. 

\bigskip
\flushpar
For every map $f\in \Cal D(X)$ we define the following functions
$$
r_n:X\to \bold{R},
$$
where $n\ge 1$ and
$$
r_n(y)=\inf \{\epsilon>0| B_{\epsilon}(y)\subset f^i(T_i) \text{ with } T_i
\text{ branch
of } f^i, i\le n, \text{ with } y\in f^i(T_i)\}.
$$
with $B_{\epsilon}(y)=(y-\epsilon,y+\epsilon)$. Furthermore define $r:X\to
\bold{R}$ to be $r=\lim r_n$.

\flushpar
The sufficient condition for the existence of acims will be formulated in
terms of the set 
$$
S=\{y\in X|r(y)>0\}.
$$

\proclaim{Theorem 4.1} Every conservative ergodic map in $\Cal D(X)$ having
$\lambda(S)>0$ exhibits
an acim.
\endproclaim

\demo{proof}
The set $S$ is backward invariant. Hence the conservativity of $f$ implies 
$|S|=1$. Let $S_{\rho}=r^{-1}(\rho,1)$ with $\rho >0$. 

\proclaim{Claim} For every compact set $I\subset S_{\rho}$ there exist finitely many pairs $\{U_i,V_i\}, i=1,...,s$, of intervals such that
\parindent=15pt
\item{1)} $I\subset \cup\{ U_i|i=1,...,s\}$;
\item{2)} for all $i=1,...,s$ $|U_i|=\rho$ and $U_i\subset V_i$ with both 
components of $V_i-U_i$
have length $\frac12 \rho$;
\item{3)} if $T$ is a branch of $f^n$ with $f^n(T)\cap (U_i\cap I)\ne 
\emptyset$ for some $i\le s$ then $V_i\subset f^n(T)$. 
\endproclaim

\demo{proof of claim}
For every $y\in I$ the interval $V_y=B_{\rho}(y)=(y-\rho,y+\rho)$ has the following property: if $T$ is a branch with $y\in f^n(T)$ then $V_y\subset f^n(T)$.

\flushpar
Consider a branch $T$ which covers a point $z\in V_y\cap I$, $z\in f^n(T)$. Because $I\subset S_\rho$ we get immediately $y\in f^n(T)$. Hence $V_y\subset 
f^n(T)$. Conclusion:
for every $y\in I$ and every branch $T$ with $f^n(T)\cap (V_y\cap I)\ne \emptyset$ we have $V_y\subset f^n(T)$.

\flushpar
Let $U_y=B_{\frac12 \rho}(y)$. By using compactness of $I$ we can cover $I$ by finitely many intervals of the form $U_y$. The corresponding pairs $\{U_y,V_y\}$ will satisfy the claim.
\enddemo

\bigskip
\flushpar
The claim implies that every compact set  $I \subset S_{\rho}$ has a natural 
finite
partition in sets $I_i=I\cap U_i$. Let us use this partitions for constructing 
$\lambda-$partitions which will allow us to apply Theorem 2.9.
   
\flushpar
As we saw $|S|=1$. Hence for every set $K$ with $|K|>0$ there exist a compact 
set $I\subset K$ of positive Lebesgue measure and a $\rho>0$ such that $I\subset S_{\rho}$. This 
observation easily implies the existence of a $\lambda-$partition 
$\Cal{G}_0$ such that for  every $I\in \Cal G_0$ $I\subset S_{\rho_I}$ for 
some  $\rho_I>0$.

\flushpar
Now partition every $I\in \Cal{G}_0$ as described above: 
$I=\cup_{i=1}^{s_I}I_i$. Define $\Cal{G}$ to be the collection consisting of 
the sets $I_i, i=1,...,s_I$ and $I\in \Cal{G}_0$.

\bigskip
\flushpar
For applying Theorem 2.9 we have to bound the distortion of the measures 
$f^n_*\lambda|I, I\in \Cal G$. 

\flushpar
Fix $I\in \Cal G$, say $I\subset S_\rho$, and let $A\subset I$. The 
definition of the sets $I\in \Cal G$ allows us to cover $f^{-n}(I)$ by 
 branches $T_1,..., T_{k_n}$, $k_n\in \bold{N}\cup\{\infty\}$, 
satisfying:   
\parindent=15pt
\item{1)} $f^n|T_i$ is diffeomorphic;
\item{2)} both components of $f^n(T_i)-\{convex-hull(I)\}$ have length 
bigger than
$\frac12 \rho$.

\flushpar
Now the Koebe-Lemma states the existence of $K>0$, only depending on $\rho$,
  such that  
$$
\frac1K\le \frac{|Df^n(x_1)|}{|Df^n(x_2)|}\le K
$$
for all $x_1,x_2\in f^{-n}(I)\cap T_i$, $i=1,...,k_n$.

\flushpar
Now
$$
\aligned
\frac{\lambda(f^{-n}(A))}{\lambda(f^{-n}(I))}&=
\frac{\sum_{i=1}^{k_n}\lambda(f^{-n}(A)\cap T_i)}{\lambda(f^{-n}(I))}\\
&=\frac{\sum_{i=1}^{k_n}\frac{\lambda(f^{-n}(A)\cap T_i)}
                             {\lambda(f^{-n}(I)\cap T_i)}
        \lambda(f^{-n}(I)\cap T_i)}
       {\lambda(f^{-n}(I))}\\
&=\frac{\sum_{i=1}^{k_n}
        \frac{\frac{\lambda(I)}
                   {\lambda(f^{-n}(I)\cap T_i)}}
             {\frac{\lambda(A)}
                   {\lambda(f^{-n}(A)\cap T_i)}}
        \frac{\lambda(A)}{\lambda(I)}\lambda(f^{-n}(I)\cap T_i)}
       {\lambda(f^{-n}(I))}.
\endaligned
$$ 
Using the mean value theorem  and the distortion result above we get
$$
\frac1K \frac{\lambda(A)}{\lambda(I)}
\le \frac{\lambda(f^{-n}(A))}{\lambda(f^{-n}(I))}
\le K \frac{\lambda(A)}{\lambda(I)},
$$
the measures $f_*^n\lambda|I$ have distortion bounded by $K$. We proved
$\lambda \in \Cal{M}(\Cal G,f)$. 
Hence Theorem 2.9 states the existence of an acim.
\qed
\enddemo

\flushpar
A map $f\in \Cal{D}(X)$ is called a {\it Markov} map if there exists a $\lambda-$partition $\Cal P$ consisting of intervals  such that
for every $I\in \Cal P$ the image $f(I)$ is a union (up to a set of measure zero) of
elements of $\Cal P$. The map $f$ is said to satisfy the Markov property with respect to $\Cal P$. Obviously these Markov maps have $r(y)>0$ for all 
$y\in \cup \Cal P$.

\proclaim{Corollary 4.2} Every conservative ergodic Markov map has an acim.
\endproclaim

\flushpar
This statement has to be compared with a theorem of Harris (see [H]) stating
the existence of infinite stationary states for certain Markov processes on
countable many state spaces. In fact examples are known of Markov process
not having a stationary state (see [D]). These examples also serve for showing
that we
cannot omit the conservativity. On the other hand  we can weaking this condition by
imposing a topological condition. By doing so we kill the metrical subtilities
and get a general existence theorem which is valid as well in the conservative
as in the dissipative case.   

\proclaim{$\sigma-$Folklore Theorem} Every finite-to-1 Markov map has an acim.
\endproclaim

\demo{proof} Suppose $f$ satisfies the Markov property with respect to $\Cal{G}_0$. Let $\Cal G$ be a $\lambda-$partition refining $\Cal{G}_0$ and consisting of intervals. These intervals $I\in \Cal{G}$ are chosen in such a way that both components of $T-I$ have length bigger than $\lambda(I)$, where $T\in \Cal{G}_0$ with $I\subset T$.

\bigskip
\flushpar
Once we proved
$$
PL_{\lambda}(\Cal{G}_0)\subset \Cal{M}(\Cal G,f)
$$
Theorem 3.1 assures the existence of an acim.

\flushpar
Fix $I\in \Cal{G}$ with $I\subset T \in \Cal{G}_0$. Furthermore let $A\subset I$. Because $f$ is a finite-to-1 Markov map the set
$f^{-n}(I)$, $n\ge o$, can be covered by finitely many intervals $T_1,..., T_{k_n}$ 
satisfying 
\parindent=15pt
\item{1)} $f^n|T_i$ is diffeomorphic;
\item{2)} $f^n(T_i)=T$;
\item{3)} there exists $T'_i\in \Cal{G}_0$ with $T_i\subset T'_i$.

\flushpar
for all $i=1,...,k_n$.
\flushpar
Take $\mu \in PL_{\lambda}(\Cal{G}_0)$. To prove that $\mu \in \Cal{M}(\Cal{G},f)$ first we have to show that $\mu(f^{-n}(I))<\infty, n\ge 0$. However this is a direct consequence of $f$ being  finite-to-1. Secondly we have to study the local distortion of the measures $f^n_*\mu$.

\flushpar
Again the Koebe-Lemma gives a constant $K>0$ such that  
$$
\frac1K\le \frac{|Df^n(x_1)|}{|Df^n(x_2)|}\le K
$$
for all $x_1,x_2\in f^{-n}(I)\cap T_i$, $i=1,...,k_n$.

\flushpar
As before
$$
\aligned
\frac{\mu(f^{-n}(A))}{\mu(f^{-n}(I))}
&=\frac{\sum_{i=1}^{k_n}
        \frac{\frac{\lambda(I)}
                   {\mu(f^{-n}(I)\cap T_i)}}
             {\frac{\lambda(A)}
                   {\mu(f^{-n}(A)\cap T_i)}}
        \frac{\lambda(A)}{\lambda(I)}\mu(f^{-n}(I)\cap T_i)}
       {\mu(f^{-n}(I))}\\
&=\frac{\sum_{i=1}^{k_n}
        \frac{\frac{\lambda(I)}
                   {\lambda(f^{-n}(I)\cap T_i)}}
             {\frac{\lambda(A)}
                   {\lambda(f^{-n}(A)\cap T_i)}}
        \frac{\lambda(A)}{\lambda(I)}\mu(f^{-n}(I)\cap T_i)}
       {\mu(f^{-n}(I))}.
\endaligned
$$ 
In the last step we used the fact that all measures in $PL_{\lambda}(\Cal{G})$ have constant densities on the elements $T\in \Cal{G}_0$ and property 3) above.
Using the mean value theorem and the distortion result above  we get
$$
\frac1K \frac{\lambda(A)}{\lambda(I)}
\le \frac{\mu(f^{-n}(A))}{\mu(f^{-n}(I))}
\le K \frac{\lambda(A)}{\lambda(I)}.
$$
We proved $\mu \in \Cal{M}(\Cal G,f)$.
\qed
\enddemo

\flushpar
The usual Folklore theorem states that every Markov map having derivative
bigger and bounded away from $1$ has an absolutely continuous invariant probability
measure. 

\flushpar
Another possible $\sigma-$Folklore theorem could be  formulated by considering
Markov maps whose branches are all mapped onto $X$.

\flushpar
As the main consequence of the $\sigma-$Folklore Theorem  we get Theorem B of
the introduction.

\proclaim{Corollary 4.3} Let $f$ be a $C^3$ map on the interval (or circle) 
satisfying
\parindent=15pt
\item{1)} $f$ has only finitely many critical points and the Schwarzian 
derivative is everywhere negative except in the critical points;
\item{2)} there exists a dense orbit;
\item{3)} the orbits of the critical points stay in a closed invariant  
set of Lebesgue
measure zero.

\flushpar
Then $f$ has a $\sigma-$finite absolutely continuous invariant measure.
\endproclaim

\demo{proof}Let $\Lambda$ be a closed invariant  set which contains the critical 
orbits. Assume it has Lebesgue measure zero. Furthermore consider the $\lambda-$partition $\Cal P$ 
consisting of the gaps of this set (the connected components of its complement). The map we are considering has
only finite critical points. Hence it is finite-to-1. In other words the map
is a finite-to-1 Markov map.
\qed      
\enddemo

\flushpar
In the unimodal case the last corollary   can be stated simpler. In 
the unimodal case the measure of the Cantor set containing the critical orbit
always has Lebesgue measure zero ([M]). So, unimodal maps which are only finitely 
renormalizable and having their critical orbit in an invariant  Cantor set have acims.

%\magnification=1200
\tolerance=3000

\bigskip
\centerline{\bf Appendix:  The Chacon-Ornstein Theorem}
\bigskip

In this section we will give a short proof of the main theorem in
conservative $\sigma-$finite ergodic theory. It is based on the classical
Birkhoff Ergodic Theorem. The proof can be summarized by: a map has an acim iff the return map on some set has an acip. 

\proclaim{Theorem A.1} Let $f:X\to X$ be  ergodic and conservative with respect to $\mu\in \Cal{B}_{\sigma}(X)$ which is a $\sigma-$finite invariant measure
for $f$. For every pair
of Riemann intergrable functions $\phi,\psi:X \to \bold{R}$ the equality
$$
\lim_{n\to \infty} \frac{\sum_{i=0}^{n-1} \phi(f^i(x))}
                        {\sum_{i=0}^{n-1} \psi(f^i(x))}
                  =\frac{\int \phi d\mu}{\int \psi d\mu}
$$
holds for $\mu-$almost every point $x \in X$. 
\endproclaim

\demo{proof}
Let $B\subset X$ having $\mu(B)<\infty$. Using the fact that $f$ is ergodic and
conservative we can write the space as a stack: $X$ is up to a set of measure
zero equal to a countable union of pairwise
disjoint
sets $B_k, k\ge 0$, where 
$$
B_k=\{x\in X|x,f(x),\dots,f^{k-1}(x)\notin B \text{ and } f^k(x)\in B\}.
$$ 
The key for the theorem is 
\demo{claim } Let $A\subset B_k$ for some $k\ge 0$. Then for $\mu-$almost every
point $x\in X$
$$
\lim_{n\to \infty}\frac{\#_n(A)}{\#_n(B)}=
                  \frac{\mu(A)}{\mu(B)}.
$$ 
Here $\#_n(U)=\#\{i=0,1,\dots,n-1| f^i(x)\in U\}$.
\enddemo

\demo{proof of claim}
The return map on $B$ is denoted by $R:B\to B$. It has the following properties
\parindent=15pt
\item{1)} $R$ is ergodic;
\item{2)} for all $x\in B$ the set $\{f(x),\dots,R(x)\}$ contains at most one
point of $A$;
\item{3)} the measure $\frac{\mu}{\mu(B)}$ is an acip for $R$.

\flushpar
Let
$$
X_A=\{x\in B|\{f(x),\dots,R(x)\}\cap A \ne \emptyset\}.
$$  
Then
\parindent=15pt
\item{4)} $\mu(X_A)=\mu(A)$.

\flushpar
The statements 1) and 2) are obvious. Let us prove 3) and 4). Consider a set
$A\subset B_k$ and define inductively the sets $A_l$ with $l\ge k$:
$$
\align
A_k&=A;\\
A_{l+1}&=f^{-1}(A_l)\cap B_{l+1}.\\
\endalign
$$
Define $R_{l+1}=f^{-1}(A_l) \cap B$ for $l\ge k$. Using induction we get
$$
\mu(A)=\sum_{l=k+1}^n \mu(R_l) +\mu(A_n)
$$
for all $n\ge k$. Applying this to $A=B$ and using the fact that almost every
point in $B$
returns to $B$ we get $\mu(B_k)\to 0$ for $k\to \infty$.

\flushpar
Consider $A\subset B_k$ and observe $X_A=\cup_{l=k+1}^\infty R_l$. Using
$\mu(A_n)\le \mu(B_n) \to 0$ we get
$$
\mu(X_A)=\mu(\bigcup_{l=1}^\infty R_l)=\sum_{l=1}^\infty \mu(R_l)=\mu(A). 
$$ 
We proved 4). Furthermore observe $R^{-1}(A)=X_A$ for every $A\subset B$.
We proved 3).

\bigskip
\flushpar
The proof of the claim is based on the Birkhoff Ergodic theorem. Consider a
point $x$ whose orbit behaves according to the invariant measure of $R$. Let
$y=f^{i_0}(x)$ be the first time when $f^i(x) \in B$. Furthermore partition the
orbit according to the returns to $B$. Then we get   
$$
\align
\lim_{n\to \infty} \frac{\#_n(A)}{\#_n(B)}&=
\lim_{n\to \infty}\frac{\#\{i=0,1,\dots,\#_n(B)-1| R^i(y)\in X_A \}}{\#_n(B)}\\
&=\frac{\mu(X_A)}{\mu(B)}=\frac{\mu(A)}{\mu(B)}.
\endalign
$$
Observe that up to time $i_0$ we hit $A$ at most $i_0$ times. Furthermore the
part of the orbit from $R^{\#_n(B)}(x)$ to $R^{\#_n(B)+1}(x)$ hit $A$ at most
once. The conservativity and ergodicity implies $\#_n(B)\to \infty$. Hence the
initial and final part of the orbit are not influencing the limit. We proved
the claim.
\enddemo

\flushpar
For general sets $A\subset X$ we get
$$
\align
\liminf_{n\to \infty}\frac{\#_n(A)}{\#_n(B)}&=
\liminf_{n\to \infty} \sum_{k\ge 0}\frac{\#_n(A\cap B_k)}{\#_n(B)}\\
&\ge \sum_{k\ge 0} \frac{\mu(A\cap B_k)}{\mu(B)}=\frac{\mu(A)}{\mu(B)}.
\endalign
$$
\enddemo

\flushpar
Using the symmetry in $A$ and $B$ we get 
$$
\frac{\mu(A)}{\mu(B)}\ge \frac1{\limsup_{n\to \infty} \frac{\#_n(B)}{\#_n(A)}}=
\liminf_{n\to \infty} \frac{\#_n(A)}{\#_n(B)} \ge \frac{\mu(A)}{\mu(B)}. 
$$
This implies, again using the symmetry, the equality
$$
\lim_{n\to \infty}\frac{\#_n(A)}{\#_n(B)}=\frac{\mu(A)}{\mu(B)}.
$$

\flushpar
The Chacon-Ornstein Theorem obviously follows for linear combinations of
indicator functions. The general statement is  a direct consequence of the
definition of Riemann integrability.
\qed
\enddemo
%\endproclaim

%\magnification=1200
\tolerance=3000

\bigskip
\centerline{\bf References}
\bigskip

\parindent=20pt
\item{[BL]} A.M.Blokh, M.Ju.Lyubich, {\it Non-existence of wandering intervals
and structure of topological attractors of one-dimensional smooth dynamics,
II.
The smooth case}, Erg. Th. \& Dyn. Sys vol. 9 (1989) 751-758. 
\item{[CE]} P.Collet, J.Eckmann, {\it Positive Liapounov exponents and
absolutely continuity for maps of the interval}, Erg. Th. \& Dyn. Sys. vol. 3
(1983) 13-46.
\item{[D]} C. Derman, {\it Contributions to the theory of denumerable Markov
chains}, Trans. Amer. 
Math. vol 79 (1955) 541-555.
\item{[F]} D.Freedman, {\it Markov Chains}, Holden-day.
\item{[GuJ]} J.Guckenheimer, S.Johnson, {\it Distortion of S-unimodal maps},
preprint (1989).
\item{[Ha]} P.R.Halmos, {\it Invariant measures}, Ann.Math. vol 48. no.3 (1947)
735-754.
\item{[H]} T.E.Harris, {\it Transient Markov chains with stationary measures},
Proc. Amer. Math. Soc. vol. 8 (1957) 937-942.
\item{[HKe1]} F.Hofbauer, G.Keller, {\it Quadratic maps without asymptotic
measure}, Comm. Math. Phys. vol 127 (1990) 319-337. 
\item{[HKe2]} F.Hofbauer, G.Keller, {\it Some remarks on recent results about S-unimodal
maps}, Ann.Inst. Henri Poincare, vol. 53 no.4 (1990) 413-425.
\item{[J]} S.Johnson, {\it Singular measures without restrictive intervals},
Comm. Math. Phys. vol 110 (1987) 185-190.
\item{[JS]} M.Jacobson, G.Swiatek, {\it Metric properties of non-renormalizable
S-unimodal maps}, preprint IHES/M/91/16.
\item{[Ka]} Y.Katznelson, {\it Sigma-finite invariant measures for smooth
mappings of the circle}, Journal d'Analyse Mat\'ematique vol. 31 (1977) 1-19.
\item{[Ke1]} G.Keller,{\it Exponents, attractors and Hopf decompositions for
interval maps}, preprint 1988.
\item{[K]} U.Krengel, {\it Ergodic Theorems}, de Gruyter.
\item{[LM]} M.Lyubich, J.Milnor, {\it The Fibonacci Unimodal map}, preprint
1991/15 SUNY, Stony Brook.
\item{[M]} M.Martens, {\it Interval Dynamics}, Thesis at Technical University
of Delft, the Netherlands.
\item{[Mi]} M.Misiurewics, {\it Absolutely continuous measures for certain maps of the interval}, Publ. Math. IHES vol. 53 (1981) 17-51.
\item{[MMS]} M.Martens, W.de Melo, S.van Strien, {\it Julia-Fatou-Sullivan
theory for real one-dimensional dynamics}, to appear in Acta Math.  
\item{[NS]} T.Nowicki, S.van Strien, {\it Invariant measures exist under a
summability condition}, preprint 1990.
\item{[O]} D.S.Ornstein, {\it On invariant measures}, Bull. AMS vol. 66 (1960)
297-300.
\item{[P]} K.Petersen, {\it Ergodic Theory}, Cambridge University Press.
\item{[S]} S.van Strien, {\it Hyperbolicity and invariant measures for general
$C^2$ intervals maps satisfying the Misiurewicz condition}, to appear in 
Comm.Math. Phys.
\item{[Sw]} G.Swiatek, {\it Bounded distortion properties of one-dimensional
maps}, preprint SUNY IMS no. 1990/10.
\bigskip
\bye